\magnification=\magstep1
\baselineskip=1.3\baselineskip

\def\wh{\widehat}
\def\eps{\varepsilon}

\def\prt{\partial}

\def\ol{\overline}
\def\R{{\bf R}}

\def\frac#1#2{{#1\over #2}}

\def\bone{{\bf 1}}

\def\sqr#1#2{{\vcenter{\vbox{\hrule height.#2pt
        \hbox{\vrule width.#2pt height#1pt \kern#1pt
           \vrule width.#2pt}
        \hrule height.#2pt}}}}

\def\qed{{\hfill $\mathchoice\sqr56\sqr56\sqr{2.1}3\sqr{1.5}3$ }}

\input epsf
\newdimen\epsfxsize

\centerline{\bf ON NODAL LINES OF NEUMANN EIGENFUNCTIONS}
\footnote{$\empty$}{\rm Research partially supported by NSF
grant DMS-0071486.}
\vskip.50truein

\centerline{\bf Krzysztof Burdzy}

\vskip.50truein

\noindent{\bf Abstract}.
\footnote{$\empty$}{\rm AMS mathematics subject classification
(2000): 35P99, 60J65}
\footnote{$\empty$}{\rm Key words or phrases: nodal line,
reflected Brownian motion, mirror coupling, eigenfunction,
Neumann problem.}
We present a new method for locating the
nodal line of the second eigenfunction
for the Neumann problem in a planar domain.
\vskip0.6truein

\noindent{\bf 1. Introduction}.
This note is inspired by recent progress on the 
``hot spots'' conjecture of J.~Rauch originally
proposed in 1974. The conjecture states
that the maximum of the second Neumann eigenfunction
in a Euclidean domain is attained on the boundary.
This formulation is somewhat vague but the counterexamples
given by Burdzy and Werner (1999) and Bass and Burdzy (2000)
show that even the weakest version of the conjecture
fails in some planar domains.
The first positive results on the ``hot spots'' conjecture
appeared in Kawohl (1985) and Ba\~nuelos and Burdzy (1999).
The last article contains a proof that the conjecture holds
in all triangles with an obtuse angle. The theorem was
generalized to ``lip domains'', i.e., planar domains between
two graphs of Lipschitz functions with the Lipschitz constant
1; this more general result was hinted at in Ba\~nuelos and Burdzy (1999)
and proved in Atar and Burdzy (2002).
See Atar (2001), Jerison and Nadirashvili (2000) 
and Pascu (2002) for some other positive results.

There is currently a significant gap between the best
positive results and the strongest counterexamples
to the ``hot spots'' conjecture. It is widely believed that
the conjecture holds in all convex domains in 
$\R^n$ for all $n\geq 2$.
But the problem remains open even for triangles
whose all angles are acute.

One of the key assumptions in some theorems
of Ba\~nuelos and Burdzy (1999), and all results
of Jerison and Nadirashvili (2000) 
and Pascu (2002), is the symmetry of the domain
with respect to a straight line. The assumption
is only needed to show that the nodal line of the
second Neumann eigenfunction lies on the line of symmetry.
Then one can replace the original Neumann problem 
on the whole domain with the mixed Dirichlet-Neumann
problem on a nodal subdomain. As it often happens
in analysis, the Dirichlet problem (in this case,
the mixed Dirichlet-Neumann problem) turns out to be
considerably easier to deal with than the original Neumann
problem. Hence, one would like to be able to find the location
of the nodal line for the second Neumann eigenfunction
in domains which are not necessarily symmetric.
If this goal is achieved, some further progress
on the ``hot spots'' conjecture may be expected.

Very few methods for finding
the nodal line are known; the following short list
is probably complete. 

\item{(a)} If a domain has a line of symmetry, the nodal
line lies on that line, subject to some natural extra assumptions.

\item{(b)} There are a handful of classes of domains, such
as rectangles and ellipses, for which explicit formulae
for eigenfunctions are known.

\item{(c)} One can find an approximate location of the
nodal line in domains which are ``long and thin'' 
in the asymptotic sense; see Jerison (2000).

\noindent
The purpose of this note is to add the following item to this list.
\medskip

\item{(d)} In some domains, one can use the probabilistic
method of ``mirror couplings'' to delineate a region
which the nodal line must intersect.

\bigskip
It should be emphasized that one cannot expect to find
the location of the nodal line of the second Neumann
eigenfunction in the sense of an explicit formula,
except in some trivial cases. What one can realistically hope to achieve
is to obtain sufficiently accurate information about the nodal line
so that this information can in turn be used 
to prove some other results of interest; 
the combination of Lemmas 2 and 3
of Burdzy and Werner (1999) is an example of such a
result on nodal lines.

Finally, we would like to mention the results of Melas (1992)
on the nodal lines in the Dirichlet case.

\bigskip
\noindent
{\bf 2. Mirror couplings and nodal lines}.
Suppose that $D$ is a planar domain (open and connected set)
with a piecewise smooth boundary.
Informally speaking, a mirror coupling is a pair $(X_t,Y_t)$
of reflecting Brownian motions in $\ol D$, such that
the line of symmetry $K_t$ between $X_t$ and $Y_t$ does not
change on any interval $(s,u)$ such that 
$X_t\notin \prt D$ and $Y_t\notin\prt D$ for all $t\in (s,u)$.
More formally,
$(X_t, Y_t)$ satisfy the following
system of stochastic differential equations,
$$
 dX=dW+dL,\quad dY=dZ+dM,\quad dZ=dW-2m\ m\cdot dW,
 \quad m=\frac{Y-X}{|Y-X|},
$$
for times less than $\zeta=\inf\{s:X_s=Y_s\}$. 
We set $Y_t = X_t$ for $t \geq \zeta$.
Here $W$ is a
planar Brownian motion,  and
$Z$ is another Brownian motion for which the increments
are mirror images of those of $W$, the mirror being the
line with respect to which
$X$ and $Y$ are symmetric. The processes $X$ and $Y$ are reflecting Brownian
motions with the singular drifts $L$ and $M$ on the boundary, i.e.,
$\int_0^\infty \bone_{ \prt D} (X_s)|dL_s| = 0$,
and similarly for $M$.
A rigorous construction of the mirror coupling
was given in Atar and Burdzy (2002) although mirror couplings
had been informally used in the past; see Atar and Burdzy (2002)
for a short review of the history of this concept.

Recall that the first eigenvalue for the Laplacian in $D$
with Neumann boundary conditions is equal to 0 and
the corresponding eigenfunction is constant. The second eigenvalue
need not be simple but its multiplicity can be only 1 or 2
(Nadirashvili (1986, 1988)).
The nodal set is the set of points in $D$ where the second
eigenfunction vanishes. If $D$ is not simply
connected, the nodal set need not be a connected
curve. The nodal set of any second Neumann eigenfunction
divides the domain $D$ into exactly two nodal
domains (connected and open sets).

\bigskip
\noindent
{\bf Theorem 1}. {\sl
Suppose that $A\subset \ol D$ is closed, $x \in D \setminus A$, and
$D_1$ is the connected component of $D\setminus A$ which contains $x$.
Assume that for some $y\in  D\setminus D_1$, the mirror $K_t$ for the coupling starting from
$(X_0,Y_0) = (x,y)$ satisfies $K_t \cap D\subset A$ for all $t\in[ 0,\zeta]$ a.s.
Then for any second Neumann eigenfunction in $D$, 
none if its nodal domains can have a closure which
is a subset of $\ol D\setminus A$ containing $x$.
}

\bigskip
\noindent
{\bf Proof}. 
Suppose that $B$ is one of the nodal domains for
a second eigenfunction and $x \in \ol B \subset \ol D \setminus A$.
Let $C_1$ be a non-empty open disc centered at a boundary
point of $B$, so small that $B_* = B\cup C_1$
satisfies $x \in \ol B_* \subset \ol D \setminus A$,
and there exists a non-empty open
disc $C_2\subset D\setminus B_*$.
Assume that for some $y\in D\setminus D_1$, if the mirror coupling starts from
$(X_0,Y_0) = (x,y)$ then $K_t \cap D \subset A$ for all $t\in[ 0,\zeta]$ a.s.
Since the stationary measure for $X_t$
is the uniform distribution in $D$, $X_t$ will
hit $C_2$ with probability one, and so it will leave
$B_*$ with probability 1. Recall that there are exactly
two nodal domains---$B$ is one of them; let $B_1$
denote the other one. Let
$\tau_X(B_*) = \inf\{t\geq 0: X_t \notin B_*\}$
and
$\tau_Y(B_1) = \inf\{t\geq 0: Y_t \notin B_1\}$.
We will argue that $\tau_X(B_*) < \tau_Y(B_1)$ a.s.
Suppose otherwise. 
By the continuity of $t\to Y_t$, $Y_{\tau_Y(B_1)} $
belongs to $\ol{\prt B_1 \setminus \prt D}$,
and so $Y_{\tau_Y(B_1)} \in \ol B\subset \ol B_*$. We have assumed
that $\tau_X(B_*) \geq  \tau_Y(B_1)$,
so $X_{\tau_Y(B_1)} \in  \ol B_*$. Since the points
$X_{\tau_Y(B_1)}$ and $Y_{\tau_Y(B_1)}$ belong
to the connected set $\ol B_*$, their line of symmetry
$K_{\tau_Y(B_1)}$ must also intersect $\ol B_*$.
We have $K_t \cap D \subset A$ for all $t\in[ 0,\zeta]$
so $A \cap \ol B_* \ne \emptyset$, which is a contradiction.
We conclude that $\tau_X(B_*) < \tau_Y(B_1)$ a.s.

Let $\mu_2 >0$ denote the second Neumann eigenvalue.
Then $\mu_2$ is the first eigenvalue for the 
mixed problem in the nodal domain $B$, with the Neumann boundary conditions
on $\prt D$ and the Dirichlet boundary conditions
on the nodal line. 
The fact that $B_*$ is strictly larger than $B$
easily implies that $\mu_2>\mu_2^*$, 
where $\mu_2^*$ is the analogous mixed
eigenvalue for $B_*$. Using the well known identification
of Brownian motion density with the heat equation solution,
we obtain from Proposition 2.1 of 
Ba\~nuelos and Burdzy (1999) that
$$\lim_{t\to \infty} P(\tau_Y(B_1) > t) e^{\mu_2 t}\in (0,\infty).$$
Since $\tau_X(B_*) < \tau_Y(B_1)$ a.s.,
$$\lim_{t\to \infty} P(\tau_X(B_*) > t) e^{\mu_2 t} < \infty,$$
and so $\mu_2^* \geq \mu_2$, but this contradicts
the fact that $\mu_2>\mu_2^*$.
Our initial assumption that
$x \in \ol B \subset \ol D \setminus A$
for some nodal domain $B$ must be false.
\qed
\bigskip

We will illustrate our main result with two examples.
Before we do that, we recall a few crucial facts about
mirror couplings from Burdzy and Kendall (2000) and 
Ba\~nuelos and Burdzy (1999).

Suppose that $D$ is a polygonal domain and 
$I$ is a line segment contained in its boundary.
Let $J$ denote the straight line containing
$I$ and recall that $K_t$ denotes the mirror line.
Let $H_t$ denote the ``hinge,'' i.e., the intersection
of $K_t$ and $J$. Note that $H_t$ need not belong
to $\prt D$. Suppose that for all $t$ in
$[t_1, t_2]$, the reflected Brownian motions
$X_t$ and $Y_t$ do not reflect from any part
of $\prt D$ except $I$. 
Let $\alpha_t$ denote the smaller of the two angles
formed by $K_t$ and $J$. Then all possible movements
of $K_t$ have to satisfy the following condition.
\bigskip

\item{({\bf M})}
The hinge does not move within the time interval
$[t_1, t_2]$, i.e., $H_t = H_{t_1}$
for all $t \in [t_1, t_2]$.
The angle $\alpha_t$ is a non-decreasing 
function of $t$ on $[t_1,t_2]$.

\bigskip
A domain $D$ with piecewise smooth boundary
can be approximated by polynomial domains $D_n$.
Mirror couplings in the approximating domains
$D_n$ converge weakly to a mirror coupling in $D$.
One can deduce which motions of the mirror
$K_t$ in $D$ are possible by analyzing all allowed
movements of mirrors in $D_n$'s and then
passing to the limit. We will not present the details
of this limit theorem here. The readers
who are concerned about full rigor should add
an extra assumption to Example 2 that $D$
is a polygonal domain.

\bigskip
\noindent
{\bf Example 1}. 
Suppose that $D$ is an obtuse triangle with 
vertices $C_1,C_2$ and $C_3$. Let $C_3$ be the vertex
with an angle greater than $\pi/2$.
Let $\ol{C_jC_k}$ denote the line segment
with endpoints $C_j$ and $C_k$, and let
$d(C_j,C_k)$ denote the distance
between these points.
Let the points $C_4, C_{10} \in \ol{C_1C_3}$, 
$C_6, C_{11} \in \ol{C_2C_3}$, and
$C_5, C_7, C_8,C_9 \in \ol{C_1C_2}$
be such that $d(C_1,C_4) = d(C_4,C_3)
= d(C_1,C_5)$,
$d(C_2,C_6) = d(C_6,C_3)
= d(C_2,C_7)$, and the following pairs of line segments
are perpendicular:
$\ol{ C_4C_8}$ and $\ol{C_2C_3}$,
$\ol{ C_6C_9}$ and $\ol{C_1C_3}$,
$\ol{ C_5C_{10}}$ and $\ol{C_1C_3}$,
$\ol{ C_7C_{11}}$ and $\ol{C_2C_3}$.
See Figure 1.

\bigskip
\vbox{
\epsfxsize=6.0in
  \qquad\hbox{\epsffile{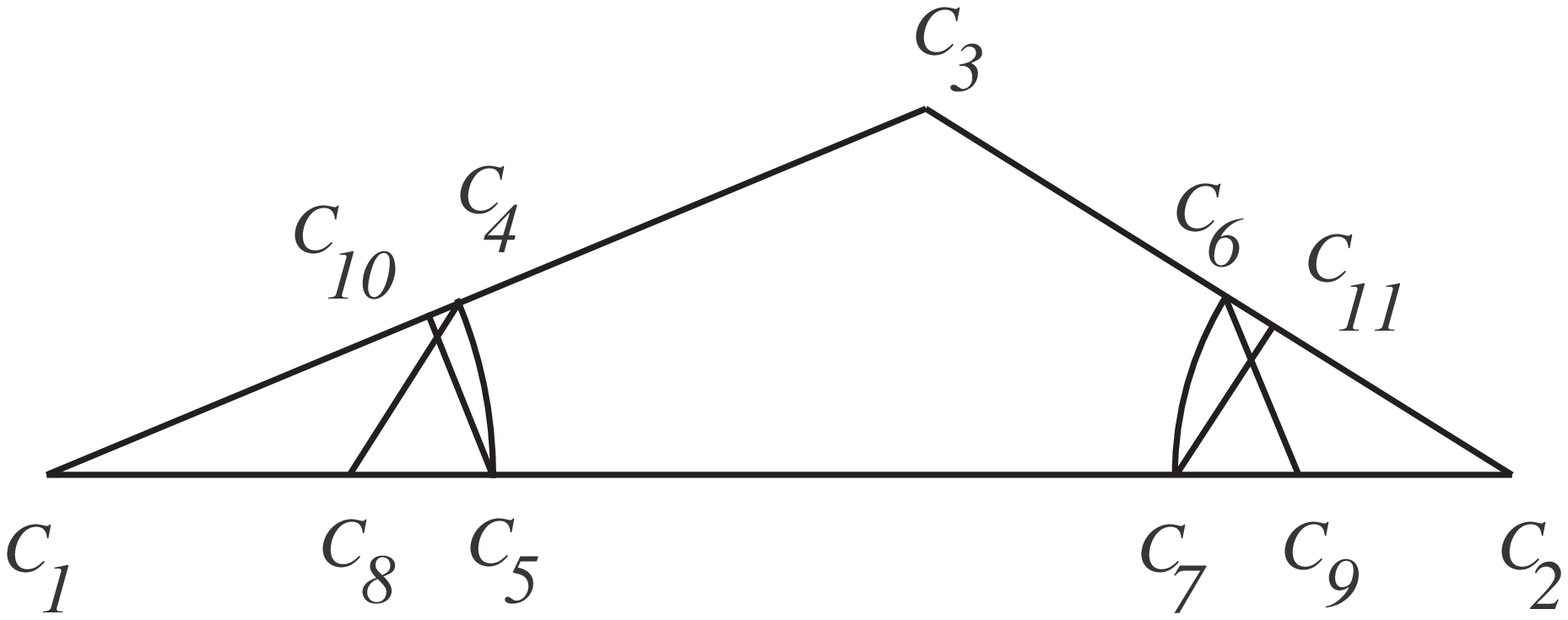}}

\centerline{Figure 1.}
}
\bigskip

Let $A$ be the closed subset of $\ol D$ whose boundary consists
of line segments $\ol{ C_3C_4}$, $\ol{C_3C_6}$ and
$\ol{ C_5C_7}$ and circular arcs $\wh{C_4C_5}$
and $\wh{C_6C_7}$. The arcs $\wh{C_4C_5}$
and $\wh{C_6C_7}$ are centered at $C_1$ and $C_2$,
respectively.
Let $A_1$ be the closure of the union of pentagons with vertices $C_3C_4C_8C_9C_6$
and $C_3C_{10}C_5C_7C_{11}$.

The second Neumann
eigenvalue is simple in obtuse triangles by a theorem of
Atar and Burdzy (2002).

\bigskip

We will show that 
\item{(i)} The nodal line for the second Neumann
eigenfunction must intersect $A$.
\item{(ii)} The nodal line lies within $A_1$.

\bigskip

Let $\alpha_t$ denote the angle formed by the mirror
$K_t$ with $\ol{ C_1C_2}$. Let $\beta_1$ and $\beta_2$
denote the angles such that if $\alpha_t = \beta_1$
then $K_t$ is perpendicular to $\ol{ C_2C_3}$
and when $\alpha_t = \beta_2$
then $K_t$ is perpendicular to $\ol{ C_1C_3}$.
An argument totally analogous to the proof of Theorem 3.1
of Ba\~nuelos and Burdzy (1999) shows that if
$\alpha_0 \in [\beta_1,\beta_2]$ then 
$\alpha_t \in [\beta_1,\beta_2]$ for all $t\geq 0$ a.s.

Consider any point $x \in D\setminus A$ to the left
of the arc $\wh{C_4C_5}$. Let $y \in \ol D$ be such
that $C_1,x$ and $y$ lie on the same line and 
$d(C_1,y) = d(C_1,C_3)$. Consider a mirror coupling
$(X_t,Y_t)$ with $X_0=x$ and $Y_0=y$. We will show
that $K_t \cap D\subset A$ for all $t\geq 0$.

We have $K_0 \cap D\subset A$, by the definition of $y$.
Assume that it is not true that $K_t \cap D\subset A$
for all $t$. We will show that this assumption
leads to a contradiction. If follows from this assumption
that for some $\eps>0$ there exists
$t$ with $d(C_1, K_t) \leq d(C_1, C_3)/2 -\eps$
or $d(C_2, K_t) \leq d(C_2, C_3)/2 -\eps$.
Let $T_1$ be the infimum of such $t$.
Without loss of generality assume that
$d(C_1, K_{T_1}) \leq d(C_1, C_3)/2 -\eps$.
Let $T_2$ be the last time $t$ before $T_1$ 
when $d(C_1, K_t) \geq d(C_1, C_3)/2 -\eps/2$.
Note that $d(C_1, K_t) \leq d(C_1, C_3)/2 -\eps/2$
for all $t\in[T_2,T_1]$.
This implies that $d(C_1, X_t) \leq d(C_1, C_3)$
and $d(C_1, Y_t) \leq d(C_1, C_3)$
for all $t\in[T_2,T_1]$. Hence, the processes
$X_t$ and $Y_t$ can reflect only on 
$\ol{ C_1C_2}$ and $\ol{ C_1C_3}$ on the time
interval $[T_2,T_1]$. According to ({\bf M}),
the distance from $C_1$ to $K_t$ cannot decrease
between times $T_2$ and $T_1$, which contradicts
the definition of these times.

Claim (i) now follows from Theorem 1.
The second claim is a consequence of the first one
and the bounds on the direction of the gradient of the
second eigenfunction proved in Theorem 3.1
of Ba\~nuelos and Burdzy (1999).

We will argue heuristically that it is impossible
to obtain much sharper estimates for the nodal line location
without imposing some extra assumptions on the obtuse triangle.
To see this, consider an obtuse triangle with two
angles very close to $\pi/2$---one slightly less
than the right angle and another one slightly larger
than that (see Figure 2). Such a triangle has a shape very close to a thin
circular sector. Let us assume that
the diameter of the triangle is equal to 1.
The nodal line in a circular sector is 
an arc with center at its vertex $C$.
The nodal line distance from $C$ is 
equal to $a_0/a_1 \approx 0.63$, where $a_0$
is the first positive zero of the Bessel function
of order 0 and $a_1$ is the first positive zero
of its derivative.
This follows from known results on eigenfunctions
in discs and estimates for Bessel function zeroes,
see Bandle (1980), p.~92. Our methods yield 0.5
as the lower bound for the distance
of the nodal line from $C$; we see that this estimate
cannot be improved beyond 0.63.

\bigskip
\vbox{
\epsfxsize=5.0in
  \quad\hbox{\epsffile{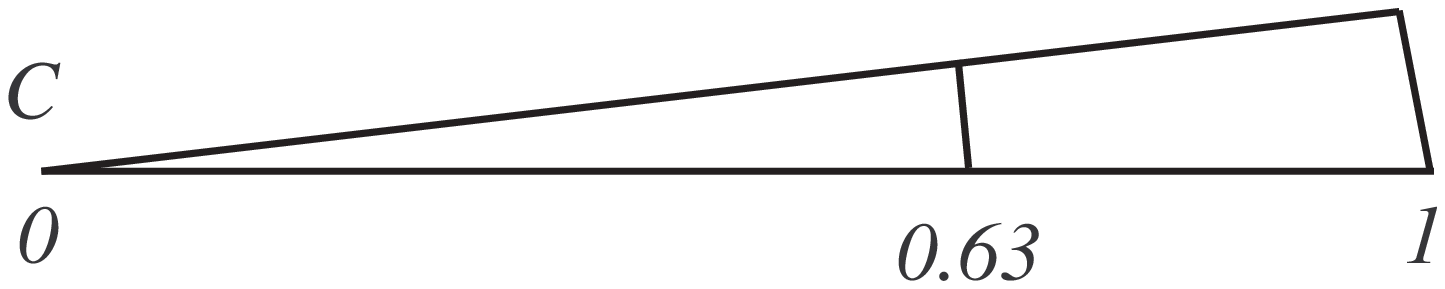}}

\centerline{Figure 2.}
}
\bigskip

We end this example with a conjecture which, if proved,
would considerably strengthen claims (i) and (ii).
\bigskip

\noindent
{\bf Conjecture 1}. {\sl
Suppose that $D$ is a triangle with vertices $C_1,C_2$ and $C_3$ and
that the angle at $C_3$ is greater than $\pi/2$.
Let $C_{12} \in \ol{C_1C_2}$ be the unique point such that
$\ol{C_{12}C_3}$ is perpendicular to $\ol{C_1C_2}$,
and let $D_1$ be the larger of the two connected
components of $D \setminus \ol{C_{12}C_3}$.
Then the nodal line for the second Neumann eigenfunction
is contained in $D_1$.
}
\bigskip

Little is known about eigenfunctions in
triangles different from equilateral; see Pinsky (1980, 1985)
for that special case. 

\bigskip
\noindent
{\bf Example 2}. 
Consider a convex domain $D$ in the plane
whose boundary contains line segments
$\{(x_1,x_2): -a \leq x_1 \leq a, x_2 = 0\}$
and 
$\{(x_1,x_2): -a \leq x_1 \leq a, x_2 = 1\}$
for some $a>0$.
Suppose that the remaining parts of $\prt D$
are contained in
$\{(x_1,x_2): -a - b \leq x_1 \leq -a, 0\leq x_2 \leq 1\}$
and
$\{(x_1,x_2): a  \leq x_1 \leq a+b, 0\leq x_2 \leq 1\}$.
We will assume that $b < a - 1/(4a)$. This means
that we have to have $a > 1/2$.
Let 
$$\eqalign{
A &= \{(x_1,x_2)\in \ol D: -b- 1/(4a)  \leq x_1 \leq b+ 1/(4a)\},\cr
A_1 &= \{(x_1,x_2)\in \ol D: -b- 1/(4a)-1  \leq x_1 \leq b+ 1/(4a)+1\}.
}$$
Figure 3 shows an example of $D$ and the corresponding
rectangle $A$.

\bigskip
\vbox{
\epsfxsize=5.0in
  \quad\hbox{\epsffile{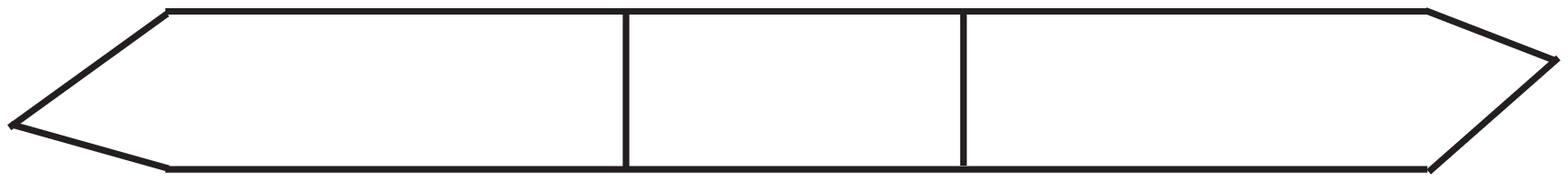}}

\centerline{Figure 3.}
}
\bigskip

We will prove that 
\item{(i)} The nodal line for any second Neumann
eigenfunction must intersect $A$.
\bigskip

We will say that $D$ is a ``lip domain'' 
if its boundary consists of graphs of two 
Lipschitz functions with the Lipschitz constant 1.
The domain depicted in Figure 3 is a lip domain.
The second eigenvalue is simple in every lip domain,
by a result of Atar and Burdzy (2002).

\bigskip

Our second claim is
\item{(ii)} If $D$ is a lip domain then
the nodal line for the Neumann eigenfunction lies within $A_1$.

\bigskip

Suppose that $x \in D\setminus A$ lies to the left of $A$.
If $x = (x_1, x_2)$ for some $x_1 \geq -a$ then let
$y=(-x_1, x_2)$. Otherwise let $y$ be any
point with the first coordinate greater than
or equal to $a$. Let
$(X_t,Y_t)$ be a mirror coupling with $X_0=x$ and $Y_0=y$. 
Our goal is to show
that $K_t \cap D\subset A$ for all $t\geq 0$.
 
It is elementary, and so it is left to the reader,
to check that for any points
$x= (x_1,x_2)$ and $y=(y_1,y_2)$ in $D$ with
$x_1 \leq -a $ and $y_1 \geq a$, the line of symmetry
for $x$ and $y$ intersects $D$ within $A$.
Hence $K_0 \cap D\subset A$, by the choice of $y$.
Next we follow the same line of thought as in Example 1.
Suppose that the condition $K_t \cap D\subset A$
fails for some $t\geq 0$. We will show that this assumption
leads to a contradiction. 
If follows from our assumption
that for some $\eps>0$ there exists
$t$ and $(z_1,z_2) \in K_t \cap \ol D$ with
$z_1 \leq - b - 1/(4a) - \eps$ or
$z_1 \geq   b + 1/(4a) + \eps$.
We can assume without loss of generality that
$z_1 \leq - b - 1/(4a) - \eps$.
Let $T_1$ be the infimum of $t$ such that
$K_t\cap \ol D$ intersects the line 
$\{(z_1,z_2): z_1= - b - 1/(4a) - \eps\}$.
Let $T_2$ be the supremum of $t\leq T_1$ 
such that
$K_t\cap \ol D$ does not intersect
$\{(z_1,z_2): z_1= - b - 1/(4a) - \eps/2\}$.
Write $X_t = (X^1_t,X^2_t)$ and $Y_t = (Y^1_t,Y^2_t)$.
The following argument applies when $\eps>0$
is sufficiently small---obviously, we can make such
an assumption on $\eps$.
Note that $X^1_t \leq a$ and $Y^1_t \leq a$
for $t\in[T_2,T_1]$. During this interval,
$Y_t$ can only reflect on the horizontal 
parts of the boundary of $D$, so according to ({\bf M}),
this reflection cannot decrease
$\inf\{z_1: (z_1,z_2) \in K_t \cap \ol D\}$.
The same holds true if $X_t$ reflects on the
horizontal parts of $\prt D$.
Suppose that $X_t$ reflects on the part
of $\prt D$ to the left of $\{(z_1,z_2): z_1 = -a\}$
which is not horizontal.
Then the hinge $H_t$ lies outside $\ol D$
and it easily follows from ({\bf M})
that both intersection points of $K_t$
with $\prt D$ move to the right. Hence,
$\inf\{z_1: (z_1,z_2) \in K_t \cap \ol D\}$
is a non-decreasing function of $t$ for 
$t\in[T_2,T_1]$. This contradicts the definition
of $T_1$ and $T_2$.

The above argument and Theorem 1 yield claim (i).
The second claim follows from the first one
and the bounds on the direction of the gradient of the
second eigenfunction given in Example 3.1
of Ba\~nuelos and Burdzy (1999).

%\bigskip
\vfill\eject
\centerline{REFERENCES}
\bigskip

\item{[1]} R. Atar (2001). Invariant wedges for a two-point
reflecting Brownian motion and the ``hot spots'' problem.
{\it Elect. J. of Probab.} 6, paper 18, 1--19.

\item{[2]} R.~Atar and K.~Burdzy (2002),
On Neumann eigenfunctions in lip domains (preprint).

\item{[3]}
C.~Bandle {\it Isoperimetric Inequalities and Applications}. 
Monographs and Studies in Mathematics, 7. Pitman, Boston, Mass.-London, 
1980.

\item{[4]} R.~Ba\~nuelos and K.~Burdzy (1999),
On the ``hot spots'' conjecture of J.~Rauch.
{\it J. Func. Anal. \bf 164}, 1--33.

\item{[5]} R.~Bass and K.~Burdzy (2000)
Fiber Brownian motion and the `hot spots' problem 
{\it Duke Math. J. \bf 105}, 25--58.

\item{[6]} K.~Burdzy and W.~Kendall (2000)
Efficient Markovian couplings: examples and counterexamples 
{\it Ann. Appl. Probab. \bf 10}, 362-409.  

\item{[7]} K.~Burdzy and W.~Werner (1999)
A counterexample to the "hot spots" conjecture.
{\it Ann. Math. \bf 149},  309--317.

\item{[8]}
D.~Jerison (2000) Locating the first nodal line in the Neumann problem. 
{\it Trans. Amer. Math. Soc. \bf 352}, 2301-2317.

\item{[9]}
D. Jerison and N. Nadirashvili (2000) The ``hot spots'' conjecture
for domains with two axes of symmetry.
{\it  J. Amer. Math. Soc. \bf 13}, 741--772.

\item{[10]} B.~Kawohl, {\it Rearrangements and Convexity of Level Sets
in PDE}, Lecture Notes in Mathematics 1150, Springer, Berlin, 1985.

\item{[11]} A.~Melas (1992) 
On the nodal line of the second eigenfunction of the Laplacian in $R\sp 2$. 
{\it J. Differential Geom. \bf 35}, 255--263. 

\item{[12]} N.S.~Nadirashvili (1986) On the multiplicity of the eigenvalues
of the Neumann problem, {\it Soviet Mathematics, Doklady, \bf 33},
281--282.

\item{[13]} N.S.~Nadirashvili (1988) Multiple eigenvalues of the Laplace operator,
{\it Mathematics of the USSR, Sbornik, \bf 133-134}, 225--238.

\item{[14]}
M. Pascu (2002)
Scaling coupling of reflecting Brownian motions
and the hot spots problem
{\it  Trans. Amer. Math. Soc.} (to appear)

\item{[15]} M.~Pinsky (1980) The eigenvalues of an equilateral triangle. 
{\it SIAM J. Math. Anal. \bf 11}, 819--827. 

\item{[16]} M.~Pinsky (1985)
Completeness of the eigenfunctions of the equilateral triangle. 
{\it SIAM J. Math. Anal. \bf 16}, 848--851.

\vskip1truecm

\baselineskip=0.7\baselineskip
\parskip=0pt

{
\obeylines
Department of Mathematics
University of Washington
Box 354350
Seattle, WA 98195-4350, USA
\medskip
Email: burdzy@math.washington.edu
http://www.math.washington.edu/\~{}burdzy/

\bye